\documentclass[a4paper]{amsart}
\usepackage{amssymb, amsmath, amsthm}

\theoremstyle{plain}
\newtheorem{thm}[subsubsection]{Theorem}
\newtheorem{cor}[subsubsection]{Corollary}
\newtheorem{prop}[subsubsection]{Proposition}
\newtheorem{lem}[subsubsection]{Lemma}

\theoremstyle{definition}
\newtheorem{deft}[subsubsection]{Definition}

\theoremstyle{remark}
\newtheorem{rem}[subsubsection]{Remark}
\newtheorem*{prob}{Problem}
\newtheorem*{exam}{Example}

\DeclareMathOperator{\Log}{Log}
\DeclareMathOperator{\Arg}{Arg}
\DeclareMathOperator{\dist}{dist}
\DeclareMathOperator{\codim}{codim}

\newcommand{\imax}{\underline{m}}

\title[Division by Flat Functions and Sectorial Extensions]{Division by Flat Ultradifferentiable Functions and Sectorial Extensions}

\author{Vincent Thilliez}

\date{}

\begin{document}

\begin{abstract}
We consider classes $ \mathcal{A}_M(S) $ of functions holomorphic in an open plane sector $ S $ and belonging to a strongly non-quasianalytic class on the closure of $ S $. In $ \mathcal{A}_M(S) $, we construct functions which are flat at the vertex of $ S $ with a sharp rate of vanishing. This allows us to obtain a Borel-Ritt type theorem for $ \mathcal{A}_M(S) $ extending previous results by Schmets and Valdivia. We also derive a division property for ideals of flat ultradifferentiable functions, in the spirit of a classical $ C^\infty $ result of Tougeron.
\end{abstract}

\maketitle

\begin{footnotesize}
\end{footnotesize}

\section*{Introduction}

Let $ S $ be an unbounded, open, plane sector with vertex at the origin. For a given sequence $ M=(M_j)_{j\geq 0} $ of positive real numbers, we consider the ultraholomorphic class of functions associated with $ M $ in $ S $, that is the class $ \mathcal{A}_M(S) $ of functions which are holomorphic in $ S $ and whose derivatives at any order $ j $ are uniformly bounded by $ C\sigma^jj!M_j $, where
$ C $ and $ \sigma $ are positive constants (depending on the function). Of course, working with such classes requires some growth and regularity assumptions on $ M $. In the present paper, we always use a standard set of assumptions, described in subsection \ref{sequences}, which ensures that a Whitney extension theorem holds for the corresponding ultradifferentiable class
$ \mathcal{C}_M(\mathbb{R}^n) $, that is the class of $ C^\infty $ functions in $ \mathbb{R}^n $ whose partial derivatives at any order $ j $ are bounded by
$ C\sigma^jj!M_j $. We shall indeed use Whitney extensions, as well as Whitney's spectral theorem, in these classes. 

At the origin, any function in $ \mathcal{A}_M(S) $ admits an asymptotic formal power series $ \sum_{j\in\mathbb{N}}\lambda_j\frac{z^j}{j!} $ with the estimate
$ \vert \lambda_j\vert\leq C\sigma^jj! M_j $ on the
coefficients. The class is said to be non-quasianalytic if it contains a non-zero function which is flat at the vertex, in other words such that $ \lambda_j=0 $ for any $ j $. The study of quasianalyticity for $ \mathcal{A}_M(S) $ has a long history, but we shall only refer, for our purpose, to the classical characterization by Korenblum \cite{Kor}, relating $ M $ and the aperture of $ S $, as recalled in subsection \ref{ultrahols}.

The present paper contains three main results. The first one, theorem \ref{sectflat}, is the key to the two others. Under the aforementioned set of assumptions on $ M $, it provides a construction of flat functions belonging to $ \mathcal{A}_M(S) $ and
admitting certain sharp estimates from above and below. In a first step, devoted to the particular case of a half-plane, the result is obtained by means of a suitable outer function, see subsection \ref{outer1}. We thank Jacques Chaumat for having suggested this
approach. In a second step, the general problem is reduced to the case of a half-plane by checking that the outer function of proposition \ref{key} behaves well
under ramification. Prior to this, a certain amount of work has to be completed in order to analyze the relationship
between $ M $ and the aperture of sectors $ S $ for which the construction is possible. We manage this by means of a growth index $ \gamma(M) $ defined and studied
in subsection \ref{growthprop}.

The second main result, theorem \ref{mythm}, deals with Borel-Ritt type theorems in the ultraholomorphic setting. Being given a formal power series $ \sum_{j\in\mathbb{N}}\lambda_j\frac{z^j}{j!} $ with 
the bounds $ \vert \lambda_j\vert\leq C\sigma^jj!M_j $, the problem consists in
finding an element of $ \mathcal{A}_M(S) $ asymptotic to that series.
Of course, the size of the sectors $ S $ in which such an extension is possible has to be related as precisely as possible to $ M $. Well-known \cite{Ra}\cite{Touglectures} for the Gevrey regularity $ \mathcal{G}^{1+\alpha} $, that is for
$ M_j=j!^\alpha $ with $ \alpha>0 $, the solution was, up to now, far from being complete in a more general setting. A noticeable step in this direction was made by Schmets and Valdivia in \cite{SV}, but the results of \cite{SV} are subject to certain limitations discussed in subsection \ref{extback}. They exclude, for instance, the case of classes which are, roughly speaking, smaller than the Gevrey class 
$ \mathcal{G}^2 $. By quite different methods, we obtain in the present work 
a ultraholomorphic Borel-Ritt theorem which, unlike the previous ones, is valid for any sequence $ M $ satisfying the aforementioned assumptions. The aperture of sectors in which the extension property holds is also sharp. The proof uses an interpolation scheme originating in classical complex analysis of several variables (see \cite{CC0} in the $ C^\infty $ case and \cite{CC1} in a special Gevrey case). The construction of theorem \ref{sectflat} plays here a crucial role and, as in \cite{SV}, the result comes with linear continuous operators. 

The last main result, theorem \ref{divflat}, deals with the factorization of flat ultradifferentiable function germs by flat factors. It is known that
any $ C^\infty $ function germ $ u $ in $ \mathbb{R}^n $ which is flat at the origin, or more generally on a germ of closed subset $ X $, can be factored as $ u=u_1u_2 $, where
both $ u_1 $ and $ u_2 $ are $ C^\infty $ and flat on $ X $. In other words, a given flat germ can always be divided by a well chosen flat germ, and the quotient is still flat. This is a particular case of a more general result of Tougeron \cite{Tougbook} recalled in
subsection \ref{divback}. Here, we address this question in the ultradifferentiable setting and we obtain the corresponding result provided $ X $ is a germ of
real-analytic submanifold.

\subsection*{Notation} For any multi-index $ J=(j_1,\ldots,j_n) $ in $ \mathbb{N}^n $, we always denote by the corresponding lower case letter $ j $ the length $ j_1+\cdots+j_n $ of $ J $. We put $ D^J =\partial^j/\partial x_1^{j_1}\cdots\partial x_n^{j_n} $.

\section{Prerequisites and basic tools}\label{tools}

\subsection{Some conditions on sequences}\label{sequences}

A sequence $ M=(M_j)_{j\geq 0} $ of real numbers is said to be
\emph{strongly regular}
if it satisfies the
following conditions, where $ A $ denotes a
positive constant:
\begin{equation}\label{Mnorm}
M_0=1 \text{ and }M\text{ is non-decreasing},
\end{equation}
\begin{equation}\label{Mlogc}
M \text{ is logarithmically convex},
\end{equation}
\begin{equation}\label{Mmodg}
M_{j+k}\leq A^{j+k}M_jM_k\text{ for any }(j,k)\in\mathbb{N}^2,
\end{equation}
\begin{equation}\label{Msnqa}
\sum_{j\geq\ell}\frac{M_j}{(j+1)M_{j+1}}\leq A\frac{M_\ell}{M_{\ell+1}}
\text{ for any }\ell\in\mathbb{N}.
\end{equation}
For any strongly regular sequence $ M $, we define its \emph{sequence of quotients}
$ m=(m_j)_{j\geq 0} $ by 
\begin{equation}\label{Mquot}
m_j=\frac{M_{j+1}}{M_j}
\text{ for any } j\in\mathbb{N}.
\end{equation}
Condition \eqref{Mlogc} amounts to saying that $ m $ 
is non-decreasing. Together with \eqref{Msnqa}, it implies
\begin{equation}\label{Multra}
\lim_{j\rightarrow\infty} m_j=\infty. 
\end{equation}
With \eqref{Mnorm}, it also implies
\begin{equation}\label{Mfast}
M_jM_k\leq M_{j+k}\text{ for any }(j,k)\in\mathbb{N}^2
\end{equation}
as well as
\begin{equation}\label{quot1}
M_j^{1/j}\leq 
m_j\text{ for any }j\in\mathbb{N}^*.
\end{equation} 
Thus, condition \eqref{Mmodg} appears as some sort of converse to \eqref{Mfast}: we refer to it as the \emph{moderate growth} condition.
It implies also a converse to \eqref{quot1}: indeed, using successively \eqref{Mlogc} and \eqref{Mmodg}, one has $ (m_j)^j\leq m_j\dots m_{2j-1}=M_{2j}/M_j\leq A^{2j} M_j $ for any integer $ j\geq 1 $, hence
\begin{equation}\label{quot2}
m_j\leq A^2 M_j^{1/j}\text{ for any } j\in\mathbb{N}^*.
\end{equation}
Finally, \eqref{Msnqa} is known as the \emph{strong non-quasianalyticity} condition;
its function-theoretical meaning
will be recalled in subsection \ref{classes}.
With $ M $ is also associated the function $ h_M $ defined on $ \mathbb{R}_+ $ by
\begin{equation*}
h_M(t)=\inf_{j\geq 0}t^j M_j\,\text{ for }t>0\,\text{ and }h_M(0)=0.
\end{equation*}
The function $ h_M $ is continuous, non-decreasing, with values in $ [0,1] $. More
precisely, by virtue of \eqref{Mnorm}, \eqref{Mlogc} and \eqref{Multra}, it is easy to see that
$ h_M(t)=t^j M_j $ for $ j\in [\frac{1}{m_j},\frac{1}{m_{j-1}}[ $ with $ j\geq 1 $, and
$ h_M(t)=1 $ for $ t\geq \frac{1}{m_0} $. In particular, $ h_M $ fully determines
$ M $ since we then have
\begin{equation*}
M_j=\sup_{t>0}t^{-j}h_M(t).
\end{equation*}
Let $ s $ be a real number, with $ s\geq 1 $. Obviously, one has 
$ \big(h_M(t)\big)^s\leq h_M(t) $. An easy but important consequence
of the moderate growth condition \eqref{Mmodg} is the
existence of a constant $ \rho(s)\geq 1 $, depending
only on $ s $ and $ M $, such that
\begin{equation}\label{hmodg}
h_M(t)\leq \big(h_M(\rho(s)t)\big)^s\text{ for any }t\in\mathbb{R}_+.
\end{equation}

\begin{exam}
A most classical example of
strongly regular sequence is given by the
Gevrey sequences $ M_j=j!^\alpha $ with $ \alpha>0 $. In this case, one
has
$ \exp(-2\alpha t^{-1/\alpha})\leq h_M(t)\leq\exp(-\frac{\alpha}{2}t^{-1/\alpha}).$
\end{exam}

\subsection{Some Carleman classes of functions}\label{classes}

Let $ M $ be a strongly regular sequence and let $ \Omega $ be an
open subset of $ \mathbb{R}^n $. For
any real $ \sigma>0 $, any complex-valued function $ f $ belonging to $ \mathcal{C}^\infty(\Omega) $ and any point $ x $ of
$ \Omega $, put
\begin{equation*}
p_\sigma(f,x)=\sup_{J\in\mathbb{N}^n}\frac{\vert D^Jf(x)\vert}{\sigma^jj!M_j}.
\end{equation*} 

\begin{lem}\label{elem}
Let $ f_1 $ and $ f_2 $ be $ C^\infty $ functions on $ \Omega $, let 
$ \sigma_1 $ and $ \sigma_2 $ be positive real numbers. Then for any point $ x $ in $ \Omega $, we have
$ p_{\sigma_1+\sigma_2}(f_1f_2,x)\leq p_{\sigma_1}(f_1,x)p_{\sigma_2}(f_2,x) $.
\end{lem}

\begin{proof} For any multi-index $ L $, we have
$ \vert D^L(f_1f_2)(x)\vert\leq p_{\sigma_1}(f_1,x)p_{\sigma_2}(f_2,x)\mathcal{S}_L $ with
\begin{equation*}
\mathcal{S}_L= \sum_{J+K=L}\frac{L!}{J!K!}\sigma_1^j\sigma_2^k j!M_j
k!M_k
\end{equation*}
by the Leibniz formula. Using the elementary estimate $ j!k!\leq (j+k)! $ and property \eqref{Mfast}, we get immediately
$ \mathcal{S}_L\leq (\sigma_1+\sigma_2)^\ell \ell!M_\ell $. The lemma follows.
\end{proof}

Now consider the space $ \mathcal{C}_{M,\sigma}(\Omega) $ of those functions $ f $ for which
$ p_\sigma(f,x) $ is uniformly bounded with respect to $ x $ in $ \Omega $, in other words, for which there exists a constant $ C_f $ such that
\begin{equation}\label{CMdeft}
\vert D^Jf(x)\vert\leq C_f\sigma^j j!M_j\text{ for any }J \in\mathbb{N}^n\text{ and any }x\in \Omega.
\end{equation}
This is a Banach space for the norm $ \Vert f\Vert_{\Omega,\sigma} $ defined as the smallest constant $ C_f $ such that \eqref{CMdeft} holds.  
We define the Carleman class $ \mathcal{C}_M(\Omega) $ as the increasing union
of all Banach spaces $ \mathcal{C}_{M,\sigma}(\Omega) $ for $ \sigma>0 $, endowed
with its natural (LB)-space topology.

\begin{rem}
In view of \eqref{CMdeft}, the sequence $ M $ conveys
the defect of analyticity of the elements of $ \mathcal{C}_M(\Omega) $.  
Depending on authors, $ M_j $ often rather denotes what appears in
\eqref{CMdeft} as $ j!M_j $. For practical reasons we prefer to separate the analytic part. Anyways, writing
$ \mathcal{M}_j=j!M_j $, it can be checked that the strong regularity of
$ M $ is equivalent to the classical set of conditions 
(M1)-(M2)-(M3) for $ \mathcal{M} $, as it appears
in \cite{BBMT} or \cite{Kom}, for instance.
\end{rem}

\begin{exam}
Taking $ M_j=j!^\alpha $ with $ \alpha>0 $, we obtain
Gevrey classes $ \mathcal{G}^{1+\alpha}(\Omega) $.
\end{exam}  

Now denote by $ \Lambda(\mathbb{N}^n) $ the space of families
$ \lambda=(\lambda_J)_{J\in\mathbb{N}^n} $ 
of complex numbers and consider  
the Borel map $ B\, :\, \mathcal{C}^\infty(\mathbb{R}^n)\longrightarrow\Lambda(\mathbb{N}^n) $,
given by $ (Bf)_J=D^Jf(0) $ for any $ f $ of $ \mathcal{C}^\infty(\mathbb{R}^n) $ and
any multi-index $ J $. 
Define $ \Lambda_{M,\sigma}(\mathbb{N}^n) $ as the space of those elements
$ \lambda $ of $ \Lambda(\mathbb{N}^n) $ for which there exists a constant
$ C_\lambda $ such that
\begin{equation}\label{GMdeft}
\vert \lambda_J\vert\leq C_\lambda\sigma^j j! M_j\text{ for any multi-index }J\in
\mathbb{N}^n.
\end{equation}
The Borel map restricts obviously as follows:
\begin{equation}\label{Taylor1}
B\, :\ \mathcal{C}_{M,\sigma}(\mathbb{R}^n)\longrightarrow \Lambda_{M,\sigma}(\mathbb{N}^n).
\end{equation}
In the same way as for functions, we define the norm $ \vert \lambda\vert_\sigma $
as the smallest possible $ C_\lambda $ in \eqref{GMdeft}. Then 
$ \Lambda_{M,\sigma}(\mathbb{N}^n) $ becomes a Banach space and the map \eqref{Taylor1} is
continuous, with norm $ 1 $. We also consider the (LB)-space $ \Lambda_M(\mathbb{N}^n) $ obtained as the inductive limit of these spaces. Thus, we have naturally 
\begin{equation}\label{Taylor2}
B\, :\ \mathcal{C}_M(\mathbb{R}^n)\longrightarrow\Lambda_M(\mathbb{N}^n).
\end{equation}
This map is surjective, as shown in particular by the following proposition, which 
summarizes results of several authors \cite{BBMT}\cite{Bruna}\cite{CC2}\cite{Pet}, not in full generality, but rather in a form adapted to our needs.

\begin{prop}\label{extclassic}
Let $ M $ be a strongly regular sequence.\\
(i) One can find a constant
$ b\geq 1 $, depending only on $ M $ and $ n $,
such that, for any real $ \sigma>0 $, there exists a linear continuous operator
\begin{equation*}
E_\sigma\, :\, \Lambda_{M,\sigma}(\mathbb{N}^n)\longrightarrow \mathcal{C}_{M,b\sigma}(\mathbb{R}^n)
\end{equation*}
satisfying $ BE_\sigma \lambda=\lambda $ for any element $ \lambda $ of $ \Lambda_{M,\sigma}(\mathbb{N}^n) $. The extensions $ E_\sigma\lambda $ can be assumed
to have compact support, contained in a prescribed neighborhood of $ 0 $.\\
(ii) For any bounded open subset $ \Omega $ of $ \mathbb{R}^n $ with 
Lipschitz boundary, one can find a constant
$ c\geq 1 $, depending on $ M $ and $ \Omega $, such that, for any real $ \sigma>0 $, there exists a linear continuous operator
\begin{equation*}
F_\sigma\, :\, \mathcal{C}_{M,\sigma}(\Omega)\longrightarrow \mathcal{C}_{M,c\sigma}(\mathbb{R}^n)
\end{equation*}
satisfying $ F_\sigma f\left\vert_\Omega\right.=f $ for any element 
$ f $ of $ \mathcal{C}_{M,\sigma}(\Omega) $.
\end{prop} 
\begin{proof}  
The proposition can be derived from sufficiently precise $ \mathcal{C}_M $ versions of 
the Whitney extension theorem, see for
instance \cite{CC2}, theorem 11 and remark 12. Part (i) corresponds to $ K=\{0\} $ and part (ii) corresponds to $ K=\overline{\Omega} $, since
the Lipschitz smoothness of $ \partial \Omega $ allows us to identify
functions in $ \mathcal{C}_M(\Omega) $ and Whitney jets of class $ \mathcal{C}_M $ on $ \overline{\Omega} $. 
\end{proof}

Since the strong regularity of $ M $ 
implies in particular the well-known Denjoy-Carleman condition of non-quasianalyticity
\begin{equation}\label{Denjoy}
\sum_{j\geq 0}\frac{M_j}{(j+1)M_{j+1}}<\infty,
\end{equation}
we know that $ \mathcal{C}_{M,\sigma}(\mathbb{R}^n) $ contains non-zero functions $ f $ which are flat at the origin, which means that $ Bf=0 $. Such a function $ f $ satisfies, for any multi-index $ K\in\mathbb{N}^n $ and any $ x\in\mathbb{R}^n $, 
\begin{equation}\label{udflat1}
\vert D^Kf(x)\vert \leq \Vert f\Vert_{\mathbb{R}^n,\sigma}(2A\sigma)^k k!M_k h_M(2A\sigma\vert x\vert).
\end{equation}
Indeed, $ D^Kf $ can be majorized by applying 
the Taylor formula at any order $ j $, together with \eqref{Mmodg} and the elementary estimate
$ (j+k)!\leq 2^{j+k}j!k! $. It suffices then to take the infimum with respect to $ j $ to
get \eqref{udflat1}.
We see in particular that any flat function $ f $ in $ \mathcal{C}_{M,\sigma}(\mathbb{R}^n) $ 
satisfies $ \vert f(x)\vert\leq \Vert f\Vert_{\mathbb{R}^n,\sigma}h_M(2A\sigma\vert x\vert) $. Using the optimal cut-off functions of Bruna \cite{Bruna}, it
is not very difficult to
construct such an $ f $ for which this estimate is sharp, in the sense
that $ \vert f(x)\vert\geq A'h_M(A'\vert x\vert) $ for some
$ A'>0 $. We shall not describe this construction, nor
the corresponding upper bounds on the successive
derivatives of $ 1/f $, since they are too crude for the
practical purposes of sections \ref{ext} and \ref{div}.
From section \ref{outer}, we shall, in appropriate circumstances,
get additional holomorphy properties, hence a much better handling of $ 1/f $.

\subsection{On growth properties of strongly regular sequences}\label{growthprop}

We study here a property of strongly regular sequences which can be viewed as
a relationship between the growth and regularity of such a sequence,
and that of suitable Gevrey sequences. 

\begin{deft}\label{Pgamma}
Let $ M $ be a strongly regular sequence, $ m $ its sequence of quotients, and let $ \gamma $ be a positive real number. We say that $ M $ \emph{satisfies property} $ (\mathsf{P}_\gamma) $ if there exist a sequence $ m'=(m'_j)_{j\geq 0} $ and a constant $ a\geq 1 $ such that 
$ (j+1)^{-\gamma}m'_j $ increases
and  $ a^{-1}m_j\leq m'_j\leq am_j $ for any $ j\in\mathbb{N} $.
\end{deft}

Notice that property $ (\mathsf{P}_\gamma) $ implies easily the estimate
\begin{equation}\label{gtgev}
a_1^j j!^{\gamma}\leq M_j\text{ for any }j\in\mathbb{N}
\end{equation}
with $ a_1=m'_0/a $, since $ M_j=m_0\cdots m_{j-1} $. 
The introduction of $ (\mathsf{P}_\gamma) $ is justified by lemma \ref{growth} below. The most important part of the lemma is a rewriting of a result of Petzsche (\cite{Pet}, corollary 1.3), for which we remind the reader about notational differences: what is denoted by $ M_j $ (resp. $ m_j $, $ m'_j $) in the present paper corresponds to $ M^*_j $ (resp. $ m^*_{j+1} $, $ n^*_{j+1} $) in \cite{Pet}. 

\begin{lem}\label{growth}
Let $ M $ be a strongly regular sequence.\\
(i) There always exists a real $ \gamma>0 $ such that property  $ (\mathsf{P}_\gamma) $ holds,\\
(ii) There exist constants $ \delta>0 $ and $ a_2>0 $ such that
\begin{equation}\label{ltgev}
M_j\leq a_2^j j!^{\delta}\text{ for any } j\in\mathbb{N}.
\end{equation}
\end{lem}
\begin{proof} Part (i) is the aforementioned result of \cite{Pet}. 
Part (ii) is obtained easily from \eqref{quot2}: choose $ \delta>0 $ such that $ A^2e^{-\delta}\leq 1. $ Then a trivial induction using \eqref{quot2} and the elementary estimate $ (\frac{j}{j+1})^{j+1}\leq e^{-1} $ yields $ M_j\leq M_1^j j^{\delta j} $. Finally we get \eqref{ltgev} by Stirling's formula.
\end{proof}

\begin{rem} In terms of Gevrey classes, \eqref{gtgev} and
\eqref{ltgev} amount to the inclusions 
$ \mathcal{G}^{1+\gamma}(\Omega)\subset \mathcal{C}_M(\Omega)\subset\mathcal{G}^{1+\delta}(\Omega) $.
\end{rem}

As a first application of lemma \ref{growth}, we gain some information on the powers
of $ M $.

\begin{lem}\label{powers}
Let $ M $ be a strongly regular sequence. Then, for any real $ s>0 $, the
sequence $ M^s=(M_j^s)_{j\geq 0} $ is also strongly regular and it satisfies
\begin{equation}\label{hpowers}
h_{M^s}(t^s)=\big(h_M(t)\big)^s\text{ for any }t\in\mathbb{R}_+.
\end{equation}
\end{lem}
\begin{proof}
The only fact which is not obvious is that $ M^s $ satisfies the strong non-quasianalyticity condition \eqref{Msnqa}. Using part (i) of lemma \ref{growth}, we have
\begin{equation*}
\begin{split}
\sum_{j\geq\ell}\frac{M_j^s}{(j+1)M_{j+1}^s} & \leq a^s\sum_{j\geq\ell}\frac{1}{(j+1)
(m'_j)^s}  = a^s\sum_{j\geq\ell}\Big(\frac{(j+1)^\gamma}
{m'_j}\Big)^s\frac{1}{(j+1)^{1+s\gamma}}\\
  & \leq a^s\Big(\frac{(\ell+1)^\gamma}
{m'_\ell}\Big)^s\sum_{j\geq\ell}\frac{1}{(j+1)^{1+s\gamma}}
\leq a^{2s}\frac{M_\ell^s}{M_{\ell+1}^s}(\ell+1)^{s\gamma}
R_\ell
\end{split}
\end{equation*}
with $ R_\ell=\sum_{j\geq\ell}(j+1)^{-(1+s\gamma)} $. The result then follows from the elementary estimate $ R_\ell\leq
(1+({s\gamma})^{-1})(\ell+1)^{-s\gamma} $. Notice that the case $ s\geq 1 $ could have been treated without lemma \ref{growth}: writing $ M_j^s/M_{j+1}^s=
(m_j)^{1-s}M_j/M_{j+1} $, it is enough to use \eqref{Msnqa} and the fact
that $ (m_j)^{1-s} $ is non-increasing in this case.
\end{proof}

We introduce now a growth index which will play a crucial role throughout the paper. 

\begin{deft}\label{defgamma}
Let $ M $ be a strongly regular sequence. We
define its \emph{growth index} $ \gamma(M) $ by
\begin{equation*}
\gamma(M)=\sup\big\{\gamma\in\mathbb{R}\, ;\, (\mathsf{P}_\gamma)\text{ holds}\big\}.
\end{equation*}
\end{deft}

By lemma \ref{growth}, this definition makes sense and we always have $ 0<\gamma(M)<\infty$. It is also easy to see that 
\begin{equation}\label{gamma}
\gamma(M^s)=s\gamma(M)\text{ for any real }s>0.
\end{equation}

\begin{exam} In the case of a Gevrey sequence $ M_j=j!^\alpha $ with $ \alpha>0 $,  one has obviously $ \gamma(M)=\alpha $. The result is the same for 
$ M_j=j!^\alpha(\Log j)^{\beta j} $ with $ \beta\in\mathbb{R} $ (either positive or negative). The case $ \alpha=1 $, $ \beta=-1 $ corresponds to the so-called
``$ 1^+$ level'' occuring in formal solutions of certain linear difference equations,
see \cite{Ecalle} (thus, $ \gamma(M)=1 $ in this case).
\end{exam}

As shown by the preceding example, property $ (\mathsf{P}_\gamma) $ is easy to test in concrete situations. However, the following lemma sheds more light
on the significance of $ \gamma(M) $. It will be used several times in the next sections.

\begin{lem}\label{light}
Let $ M $ be a strongly regular sequence. Then for any real number $ \gamma $ with
$ 0<\gamma<\gamma(M) $, there exist a constant $ a\geq 1 $ and a strongly regular sequence $ M' $ such that 
$ (j!^{-\gamma}M'_j)_{j\geq 0} $ is strongly regular and 
$ a^{-j}M_j\leq M'_j\leq a^jM_j $ for any $ j\in\mathbb{N} $.
\end{lem}
\begin{proof}
Choose a real $ \delta $ with
$ \gamma<\delta<\gamma(M) $. By definition of $ \gamma(M) $, property
$ (\mathsf{P}_\delta) $ holds, hence
one can find a constant $ a\geq 1 $ 
and a sequence $ m'=(m'_j)_{j\geq 0} $ such that $ (j+1)^{-\delta}m'_j $ increases and
$ a^{-1}m_j\leq m'_j\leq am_j $ for any $ j\in\mathbb{N} $. 
Clearly we can also assume $ m'_0\geq 1 $. 
Put $ M'_0=1 $ and $ M'_j=m'_0\cdots m'_{j-1} $ for $ j\geq 1 $. By straightforward verifications, $ M' $ is strongly regular and it satisfies the estimate
$ a^{-j}M_j\leq M'_j\leq a^j M_j $ for any $ j\in\mathbb{N} $. 
Next we show that the
sequence $ M'' $ given by $ M''_j=j!^{-\gamma}M'_j $ is strongly regular. It 
is easy to see that it satisfies conditions \eqref{Mnorm} and \eqref{Mlogc}. The moderate growth
property \eqref{Mmodg} for $ M'' $ is a consequence of the moderate growth property
of $ M' $ and of the elementary estimate $ j!k!\leq (j+k)!\leq 2^{j+k}j!k! $.
Finally, for any integer
$ \ell\geq 0 $, we have
\begin{equation*}
\begin{split}
\sum_{j\geq\ell}\frac{M''_j}{(j+1)M''_{j+1}}
& =\sum_{j\geq\ell}\frac{1}{(j+1)^{1-\gamma}m'_j}=\sum_{j\geq\ell}\frac{(j+1)^\delta}
{m'_j}\frac{1}{(j+1)^{1+\delta-\gamma}}
 \\
& \leq \frac{(\ell+1)^\delta}
{m'_\ell}\sum_{j\geq\ell}\frac{1}{(j+1)^{1+\delta-\gamma}}
= \frac{M''_\ell}{M''_{\ell+1}}(\ell+1)^{\delta-\gamma}R_\ell
\end{split}
\end{equation*}
with $ R_\ell=\sum_{j\geq\ell}(j+1)^{-(1+\delta-\gamma)} $.
The strong non-quasianalyticity condition \eqref{Msnqa} for $ M'' $ follows, since 
we have $ R_\ell\leq 
(1+({\delta-\gamma})^{-1})(\ell+1)^{\gamma-\delta} $.
\end{proof}

\section{Sectorial flatness}\label{outer}

In this section, we construct holomorphic functions in plane sectors with precise flatness properties at the vertex. For any real $ \gamma>0 $, put 
\begin{equation*}
S_\gamma
=\Big\{ z\in\Sigma\ ;\ \vert \Arg z\vert<\gamma\frac{\pi}{2}\Big\},
\end{equation*}
where $ \Sigma $ denotes the
Riemann surface of the logarithm and
$ \Arg $ the principal determination of the argument. We consider first the special case
of the right half-plane $ S_1=\{z\in\mathbb{C}\, ;\, \Re z>0\} $.

\subsection{A construction of outer functions}\label{outer1}
We begin by stating two auxiliary technical lemmas. The
second one is a familiar logarithmic integral condition: this will be, in fact, the starting point of our construction.

\begin{lem}\label{weight}
Let $ N $ be a strongly regular sequence with $ \gamma(N)>1 $. Then there exist real constants $ b_1 $ and $ b_2 $, with $ b_1>0 $, such
that, for any real $ u>0 $,
\begin{equation*}
\int_0^1\Log h_N(su)ds\geq \Log h_N(b_1u)+b_2.
\end{equation*}
\end{lem}
\begin{proof} Let $ N' $ be the sequence associated with $ N $ and with $ \gamma=1 $ by lemma \ref{light}, so that $ (j!^{-1}N'_j)_{j\geq 0} $ is strongly regular and
\begin{equation}\label{divgev1}
a^{-j}N_j\leq N'_j\leq a^jN_j \textrm{ for any } j\in \mathbb{N}.
\end{equation}
Define now
\begin{equation*}
\omega(t)=\sup_{j\geq 0}\Log\Big(\frac{t^j}{N'_j}\Big)\text{ for }t>0,\quad\omega(0)=0.
\end{equation*}
Since $ (j!^{-1}N'_j)_{j\geq 0} $ satisfies the strong non-quasianalyticity condition \eqref{Msnqa}, a result
of Komatsu (\cite{Kom}, proposition 4.4, equation 4.14) shows that one can find 
a constant $ b_3>0 $ such that
\begin{equation}\label{omega}
\int_1^{+\infty}\frac{\omega(ty)}{t^2}dt\leq b_3\omega(y)+b_3\text{ for any }y>0
\end{equation}
(in fact, the strong regularity of $ (j!^{-1}N'_j)_{j\geq 0} $ implies that $ \omega $ is a strong weight function in the sense of \cite{BBMT}).
Besides \eqref{omega}, it is clear by \eqref{divgev1} that one has
\begin{equation*}
-\Log\Big(h_N\Big(\frac{a}{t}\Big)\Big)\leq\omega(t)\leq -\Log\Big(h_N\Big(\frac{1}{at}\Big)\Big).
\end{equation*}
Putting $ t=1/s $ and $ y=1/u $ in \eqref{omega} then yields easily the result.
\end{proof}

\begin{lem}\label{logint}
Let $ N $ be a strongly regular sequence with $ \gamma(N)>1. $ Then we have
\begin{equation*}
\int_{-\infty}^{+\infty}\frac{\Log h_N(\vert t\vert)}{1+t^2}dt>-\infty.
\end{equation*}
\end{lem}
\begin{proof}
Pick a real $ \gamma $ with $ 1<\gamma<\gamma(N) $. By definition of $ \gamma(N) $, and
by property \eqref{gtgev} for the sequence $ N $, we have 
$ a_1^j j!^\gamma\leq N_j $ for any $ j\in\mathbb{N} $. Multiplying this inequality by $ \vert t\vert^j $ and taking the infimum with respect to $ j $, we derive $ \exp(-b_4\vert t\vert^{-1/\gamma})\leq h_N(\vert t\vert) $ for some suitable $ b_4>0 $. The lemma follows, since $ 1/\gamma<1 $.
\end{proof}

We are now ready to work out the key construction.

\begin{lem}\label{key}
Let $ N $ be a strongly regular sequence with $ \gamma(N)>1 $. Then
there exist a function $ F $ holomorphic in the right half-plane $ S_1 $ and
constants $ b_5>0 $, $ b_6>0 $ and $ b_7>0 $, depending only on $ N $, such
that 
\begin{equation}\label{Fbounds}
b_5h_N(b_6\Re w)\leq\vert F(w)\vert \leq h_N(b_7\vert w\vert)\text{ for any }w\in S_1.
\end{equation}
\end{lem}
\begin{proof} We use the classical construction of outer functions in $ H^p $ spaces,
see e.g. \cite{Duren}, theorem 11.6, or \cite{Hof}, chapter 8. For $ w\in S_1 $, we put 
\begin{equation*}
F(w)=\exp\left(\frac{1}{\pi}\int_{-\infty}^{+\infty}\Log h_N(\vert t\vert)
\frac{itw-1}{it-w}\frac{dt}{1+t^2}\right).
\end{equation*}
Thanks to lemma \ref{logint} and to the boundedness of $ h_N $, it is known that
$ F $ is a bounded holomorphic function in $ S_1 $. 
For $ w\in S_1 $, put now $ w=u+iv $ (thus, $ u>0 $). Then we have
\begin{equation}\label{up1}
\vert F(w)\vert =\exp\left(\frac{1}{\pi}\int_{-\infty}^{+\infty}\Log h_N(\vert t\vert)
\frac{u}{(t-v)^2+u^2}dt\right).
\end{equation}
Remark that $ \int_{\vert t-v\vert\leq u}\frac{u}{(t-v)^2+u^2}dt=\frac{\pi}{2} $. Since $ \vert t-v\vert\leq u $ implies $ \vert t\vert\leq 2\vert w\vert $, we derive
\begin{equation}\label{up2}
\int_{\vert t-v\vert\leq u}\Log h_N(\vert t\vert)\frac{u}{(t-v)^2+u^2}dt \leq
\frac{\pi}{2}\Log h_N(2\vert w\vert).  
\end{equation}
Since $ h_N(\vert t\vert)\leq 1 $ for any $ t $, we also have
\begin{equation}\label{up3}
\int_{\vert t-v\vert> u}\Log h_N(\vert t\vert)\frac{u}{(t-v)^2+u^2}dt\leq 0
\end{equation}
By \eqref{up1}, \eqref{up2} and
\eqref{up3}, we obtain $ \vert F(w)\vert\leq \sqrt{h_N(2\vert w\vert)} $. 
Using property \eqref{hmodg} with $ M=N $, we get
the upper bound in \eqref{Fbounds}, with $ b_7=2\rho(2) $. The proof of the
lower bound goes as follows. Putting $ t=us $ in \eqref{up1} gives
\begin{equation}\label{low1}
\vert F(w)\vert =\exp\left(\frac{1}{\pi}\int_{-\infty}^{+\infty}\Log h_N(u\vert s\vert)
\frac{u^2}{(us-v)^2+u^2}ds\right).
\end{equation}
For $ \vert s\vert\geq 1 $, one has $ \Log h_N(u\vert s\vert)\geq\Log h_N(u) $, hence \begin{equation*}
\int_{\vert s\vert\geq 1}\Log h_N(u\vert s\vert)\frac{u^2}{(us-v)^2+u^2}ds\geq\Log h_N(u)  \int_{\vert s\vert\geq 1}\frac{u^2}{(us-v)^2+u^2}ds.
\end{equation*}
Since $ \Log h_N(u)\leq 0 $ and
$ \int_{\vert s\vert\geq 1}\frac{u^2}{(us-v)^2+u^2}ds\leq\pi $, this implies
\begin{equation}\label{low2}
\int_{\vert s\vert\geq 1}\Log h_N(u\vert s\vert)\frac{u^2}{(us-v)^2+u^2}ds\geq
\pi\Log h_N(u).
\end{equation}
One has also $ \frac{u^2}{(us-v)^2+u^2}\leq 1 $ and 
$ \Log h_N(u\vert s\vert)\leq 0 $ for any $ u $, $ v $ and $ s $. Invoking lemma \ref{weight}, we get therefore
\begin{equation}\label{low3}
\begin{split}
\int_{\vert s\vert\leq 1}\Log h_N(u\vert s\vert)\frac{u^2}{(us-v)^2+u^2}ds & \geq
\int_{\vert s\vert\leq 1}\Log h_N(u\vert s\vert)ds\\
& \geq 2(\Log h_N(b_1u)+b_2).
\end{split}
\end{equation}
From \eqref{low1}, \eqref{low2} and \eqref{low3}, we derive easily
$ \vert F(w)\vert\geq b_5(h_N(b_8u))^2 $ with $ b_5=\exp(2b_2/\pi) $ and
$ b_8=\min(1,b_1) $. By virtue of property \eqref{hmodg} for $ N $, this yields
the desired lower bound, with $ b_6=b_8/\rho(2) $.
\end{proof}

\subsection{Background on ultraholomorphic functions}\label{ultrahols}

For any open subset $ \Omega $ of $ \mathbb{C}^n $, we use the
standard identification between $ \mathbb{C}^n $ and $ \mathbb{R}^{2n} $ to
consider the spaces $ \mathcal{C}_{M,\sigma}(\Omega) $ and $ \mathcal{C}_M(\Omega) $ as introduced
in subsection \ref{classes}. We put  
$ \mathcal{A}_{M,\sigma}(\Omega)=\mathcal{H}(\Omega)\cap \mathcal{C}_{M,\sigma}(\Omega) $, where
$ \mathcal{H}(\Omega) $ denotes the space of holomorphic functions in $ \Omega $. 
Clearly, 
$ \mathcal{A}_{M,\sigma}(\Omega) $
is a closed subspace of $ \mathcal{C}_M(\Omega) $. In the same way, we put
$ \mathcal{A}_M(\Omega)=\mathcal{H}(\Omega)\cap \mathcal{C}_M(\Omega) $. 

Now consider the sectors $ S_\gamma $ introduced at the beginning of the current section.
For $ \gamma<2 $, one can use the preceding definition with $ \Omega=S_\gamma $ since
in this case $ S_\gamma $ is an open subset of $ \mathbb{C} $. It is then easy to see that any element $ f $ of $ \mathcal{A}_M(S_\gamma) $
extends continuously, together with all its derivatives, to the closure
of $ S_\gamma $. In particular, $ f $ has a Taylor series at $ 0 $.
For $ \gamma\geq 2 $, in other words for sectors on the Riemann surface $ \Sigma $, one can  similarly define $ \mathcal{A}_M(S_\gamma) $ 
as the space of holomorphic functions $ f $ in $ S_\gamma $ whose derivatives at any order
$ j $ are uniformly bounded by $ C_f\sigma^jj!M_j $ for suitable constants $ C_f $ and
$ \sigma $. The Taylor series still makes sense, since all the restrictions
of $ f $ to subsectors of aperture smaller than $ 2\pi $ have the same
expansion at $ 0 $. In all cases, one has thus a Borel map 
(abusively still denoted by $ B $)
\begin{equation*}
\begin{split}
B\, :\ & \mathcal{A}_M(S_\gamma) \longrightarrow \Lambda_M(\mathbb{N}).\\
 &  \quad f\longmapsto (f^{(j)}(0))_{j\in\mathbb{N}}
\end{split}
\end{equation*}
As for $ \mathcal{C}_M(\mathbb{R}^n) $, an element $ f $ of
$ \mathcal{A}_M(S_\gamma) $ is said to be flat if $ Bf=0 $, and the class
$ \mathcal{A}_M(S_\gamma) $ is said to be non-quasianalytic if it contains a non-zero
function $ f $ which is flat at $ 0 $. A classical work of
Korenblum provides a necessary and sufficient condition for this property (\cite{Kor}, theorem 3 and remark 1). 
In our setting, thanks to \eqref{quot1}, \eqref{quot2} and Stirling's formula, Korenblum's result can be stated as follows: $ \mathcal{A}_M(S_\gamma) $ is non-quasianalytic if and only if  
\begin{equation}\label{Koren}
\sum_{j\geq 1}\left(\frac{M_j}{(j+1)M_{j+1}}\right)^{\frac{1}{\gamma+1}}<\infty.
\end{equation} 

Just as the Denjoy-Carleman condition appears as a limit case of \eqref{Koren} (putting $ \gamma=0 $), it turns out that the strong non-quasianalyticity condition \eqref{Msnqa} appears as the limit case of another condition \eqref{sKoren}, stronger than \eqref{Koren}, and related to the growth index $ \gamma(M) $  by means of the following lemma. 
 
\begin{lem}\label{classqa}
Let $ M $ be a strongly regular sequence and let $ \gamma $ be a real number with
$ 0<\gamma<\gamma(M) $. Then there exists a constant $ b_9>0 $ such that 
\begin{equation}\label{sKoren}
\sum_{j\geq\ell}\left(\frac{M_j}{(j+1)M_{j+1}}\right)^{\frac{1}{\gamma+1}}\leq
b_9(\ell+1)\left(\frac{M_\ell}{(\ell+1)M_{\ell+1}}\right)^{\frac{1}{\gamma+1}}
\text{ for any }\ell\in\mathbb{N}.
\end{equation} 
\end{lem}
\begin{proof}
Consider the sequence $ M' $ associated with $ M $ and $ \gamma $ by lemma \ref{light}. Since we have $ a^{-1}m'_j\leq m_j\leq am'_j $ for some constant $ a\geq 1 $, it suffices to prove \eqref{sKoren} with $ M $ replaced by $ M' $. Put $ M''_j=j!^{-\gamma}M'_j $ and
$ M'''_j=(M''_j)^{\frac{1}{\gamma+1}} $. Since $ M'' $ is strongly regular, lemma \ref{powers} shows that $ M''' $ is also strongly regular. Writing the strong non-quasianalyticity property \eqref{Msnqa} for the sequence $ M''' $, we obtain precisely the desired estimate.
\end{proof}

\begin{rem} 
If we replace $ \gamma $ by a natural integer $ r $ in \eqref{sKoren}, we recover the so-called property $ (\gamma_{r+1}) $ of Schmets and Valdivia \cite{SV}.
\end{rem}

\subsection{A result on sectorially flat functions}\label{sharpsectflat}
Lemma \ref{classqa} shows in particular that
$ \mathcal{A}_M(S_\gamma) $ is non-quasianalytic provided $ \gamma<\gamma(M) $ (this is no longer true for $ \gamma\geq\gamma(M) $, as the Gevrey case shows). The construction of lemma \ref{key} allows us to state a much more precise result, since we can now obtain flat functions with sharp estimates, as announced above. 

\begin{thm}\label{sectflat} 
Let $ M $ be a strongly regular sequence and let $ \gamma $ be a real number, with
$ 0<\gamma<\gamma(M) $. There exists a function $ G $ belonging to
$ \mathcal{A}_M(S_\gamma) $ such that, for any $ z\in S_\gamma $, we have the estimate
\begin{equation}\label{sharp}
\kappa_1h_M(\kappa_2\vert z\vert)\leq \vert G(z)\vert\leq h_M(\kappa_3\vert z\vert)
\end{equation}
where
$ \kappa_1 $, $ \kappa_2 $ and $ \kappa_3 $ are positive constants depending
only on $ M $ and $ \gamma $.
\end{thm}
\begin{proof}
Pick two real numbers $ \delta $ and $ s $ with $ \gamma<\delta<\gamma(M) $ and
$ s\delta<1<s\gamma(M) $. Consider
the sequence $ N=M^s $. We know from lemma \ref{powers} that $ N $ is strongly
regular, and by \eqref{gamma} we have $ \gamma(N)=s\gamma(M)>1 $. It is
thus possible to apply lemma \ref{key}. Consider the function $ F $ provided by
the lemma and put
\begin{equation*}
G(z)=F(z^s)\text{ for }z\in S_\delta.
\end{equation*}
This makes sense since $ z\longmapsto w=z^s $ maps holomorphically 
$ S_\delta $ into the subsector $ S_{s\delta} $ of $ S_1 $. We shall show that the restriction of $ G $ to the
subsector $ S_\gamma $ of $ S_\delta $ has all the desired properties. We claim first that the
estimates \eqref{sharp} hold for $ z $ in $ S_\delta $ (hence in
$ S_\gamma $). The lower estimate is obtained from \eqref{Fbounds} by the following arguments: for $ w\in S_{s\delta} $, one has 
$ \Re w\geq b_{10}\vert w\vert $, with $ b_{10}=\cos(s\delta\frac{\pi}{2})>0 $. For
$ z\in S_\delta $, one has thus $ h_N(b_6\Re z^s)\geq h_N((b_{11}\vert z\vert)^s) $ with $ b_{11}=(b_6b_{10})^{1/s} $. One has also
$ h_N((b_{11}\vert z\vert)^s)=
(h_M(b_{11}\vert z\vert))^s $, as observed in \eqref{hpowers}. For $ s\leq 1 $, we
derive immediately $ h_N(b_6\Re z^s)\geq h_M(b_{11}\vert z\vert) $. For $ s\geq 1 $, we
use \eqref{hmodg} to obtain the same estimate, where the value of $ b_{11} $ is
divided by $ \rho(s) $. In any case, we get the desired lower bound in \eqref{sharp}, with
$ \kappa_1=b_5 $ and $ \kappa_2=b_{11} $. The proof of the upper estimate goes along the same lines and we
skip the details. Finally, we have to show that $ G $ belongs to
$ \mathcal{A}_M(S_\gamma) $. Choose a real $ \varepsilon $ with
$ 0<\varepsilon<\min(1,\delta-\gamma)\frac{\pi}{2} $. Then for any
$ z\in S_\gamma $, the closed disc of center $ z $ and radius
$ (\sin\varepsilon)\vert z\vert $ lies in $ S_\delta $. Since \eqref{sharp} has been shown to hold on $ S_\delta $, the Cauchy formula yields easily
\begin{equation}\label{cauchy1}
\vert G^{(j)}(z)\vert\leq \frac{j!}{((\sin\varepsilon)\vert z\vert)^j}
h_M\big(\kappa_3(1+\sin\varepsilon)\vert z\vert\big)\text{ for any }j\in\mathbb{N}.
\end{equation}
Since $ h_M(t)\leq t^j M_j $ for any $ j $, we derive
$ \vert G^{(j)}(z)\vert\leq b_{12}^j j!M_j $ for any $ j $, with
$ b_{12}=\kappa_3\big(1+(\sin\varepsilon)^{-1}\big) $. This completes the proof.
\end{proof}

The following supplement to theorem \ref{sectflat} will be useful in section \ref{ext}.

\begin{lem}\label{flatuplo}
The function $ G $ of theorem \ref{sectflat} satisfies, for any $ z\in S_\gamma $ and
any $ j\in\mathbb{N} $, 
\begin{equation}\label{flatup}
\vert G^{(j)}(z)\vert \leq b_{13}^j j!M_j\, h_M(b_{14}\vert z\vert)
\end{equation}
\begin{equation}\label{flatlo}
\Big\vert\Big(\frac{1}{G}\Big)^{(j)}(z)\Big\vert \leq b_{15}b_{16}^j j!M_j
\big(h_M(b_{17}\vert z\vert)\big)^{-1}
\end{equation}
with positive constants $ b_{13} $ to $ b_{17} $ depending only on $ M $ and $ \gamma $.
\end{lem}
\begin{proof} Going back to \eqref{cauchy1}, and using \eqref{hmodg} with $ s=2 $, we
obtain easily \eqref{flatup} with $ b_{14}=\rho(2)\kappa_3(1+\sin\varepsilon) $ and
$ b_{13}=b_{14}(\sin\varepsilon)^{-1} $. The proof of \eqref{flatlo} follows the same pattern: using the lower bound of $ G $ in $ S_\delta $ and the Cauchy formula on the same disc as before, we get
\begin{equation*}
\Big\vert\Big(\frac{1}{G}\Big)^{(j)}(z)\Big\vert 
\leq \kappa_1^{-1}\frac{j!}{((\sin\varepsilon)\vert z\vert)^j}
\big(h_M(b_{18}\vert z\vert)\big)^{-1}
\end{equation*}
with $ b_{18}=\kappa_2(1-\sin\varepsilon) $. Using once more \eqref{hmodg}, we remark that
\begin{equation*}
\big(h_M(b_{18}\vert z\vert)\big)^{-1}= \frac{h_M(b_{18}\vert z\vert)}{\big(h_M(b_{18}\vert z\vert)\big)^2}\leq 
\frac{h_M(b_{18}\vert z\vert)}{h_M(b_{17}\vert z\vert)}\leq \frac{(b_{18}\vert z\vert)^j M_j}{h_M(b_{17}\vert z\vert)}
\end{equation*}
with $ b_{17}=b_{18}/\rho(2) $. Hence we
obtain the result with $ b_{15}=1/\kappa_1 $ and $ b_{16}=b_{18}(\sin\varepsilon)^{-1} $. 
\end{proof}

\begin{rem}\label{easier}
Theorem \ref{sectflat} is much easier for Gevrey classes. Indeed, for $ M_\ell=\ell!^\alpha $ with
$ \alpha>0 $, it is
possible to check directly that $ G(z)=\exp(-z^{-1/\alpha}) $ has all the required
properties. 
The key fact is that, in this Gevrey
setting, $ h_M(t) $ is comparable (up to scaling constants) to $ \exp(-t^{-1/\alpha}) $, which appears directly as the restriction of $ G $ to $ \mathbb{R}_+ $. There is no such
explicit estimate for general sequences $ M $, and the preceding work amounts to constructing a function which plays a similar role for the corresponding $ h_M $.
\end{rem}

\section{Sectorial extensions}\label{ext}

\subsection{Background on Borel-Ritt type theorems}\label{extback}

The well-known Borel-Ritt theorem states that for any element $ \lambda $ of
$ \Lambda(\mathbb{N}) $ and any sector $ S_\gamma $ with $ 0<\gamma<2 $, one can find
a holomorphic function on $ S_\gamma $ having $ \sum_{j\in\mathbb{N}}\lambda_j\frac{z^j}{j!} $ as asymptotic expansion at $ 0 $. It
implies the Borel theorem in its most classical form, that is the surjectivity of the map $ B\, :\, \mathcal{C}^\infty(\mathbb{R})\longrightarrow\Lambda(\mathbb{N}) $. 
Since the Borel theorem admits ultradifferentiable versions,
it is natural to ask whether the Borel-Ritt theorem has ultraholomorphic analogues. 
In our context, the problem can be stated as follows:

\begin{prob}
Find conditions, relating the strongly regular sequence
$ M $ and the real number $ \gamma $, which ensure the surjectivity of $ B\, :\, \mathcal{A}_M(S_\gamma)\longrightarrow \Lambda_M(\mathbb{N}) $. 
\end{prob}

In the typical 
Gevrey case, the answer to this question
is well-known as a basic tool in the asymptotic theory of differential equations: see e.g. \cite{Ra} or \cite{Touglectures} and
the references therein. Precisely, when $ M_j=j!^\alpha $, the map $ B\, :\, \mathcal{A}_M(S_\gamma)\longrightarrow \Lambda_M(\mathbb{N}) $ is surjective if and only if $ \gamma<\alpha $, that is
$ \gamma<\gamma(M) $. The classical proof is based on explicit constructions which are specific to Gevrey classes, as those mentioned in remark \ref{easier}. For more general
sequences $ M $, despite the quite particular results of \cite{BrakImm}, very few things were known until the recent article \cite{SV} of Schmets and Valdivia. Theorem 5.8 of
\cite{SV} is more particularly related to our problem; we recall it briefly for the reader's convenience. Being given an integer $ r\geq 0 $, denote by 
$ \mathcal{D}_{r+1} $ the space of $ \mathcal{C}^\infty $ functions $ f $ on the real line, supported in $ [-1,1] $, for which one can find constants $ c_1 $ and $ c_2 $ such that, for any $ j\in\mathbb{N} $, one has
$ \sup_{x\in\mathbb{R}}\vert f^{((r+1)j)}(x)\vert\leq c_1c_2^jj!M_j $ and $ f^{((r+1)j+k)}(0)=0 $ for $ k=1,\ldots,r $. Then a sufficient condition for the Borel map $ B\, :\, \mathcal{A}_M(S_\gamma)\longrightarrow\Lambda_M(\mathbb{N}) $ to be surjective for any real number $ \gamma $ with $ 0<\gamma<r $, is that the map
$ R\, :\, \mathcal{D}_{r+1}\longrightarrow\Lambda_M(\mathbb{N}) $ defined by
$ Rf=(f^{((r+1)j)}(0))_{j\in\mathbb{N}} $ be itself surjective. 
Moreover, the extensions are given by linear continous operators between Banach spaces
$ \Lambda_{M,\sigma}(\mathbb{N}) $ and $ \mathcal{A}_{M,d\sigma}(S_\gamma) $ for some suitable constant $ d\geq 1 $, depending only on $ M $ and $ \gamma $. As pointed out in \cite{SV}, the preceding sufficient condition implies in particular \eqref{sKoren} with $ \gamma=r $. It is thus generally false for $ r\geq\gamma(M) $. Conversely, using lemmas \ref{light}, \ref{classqa} and proposition \ref{extclassic}, it can be shown that the condition holds provided
$ r<\gamma(M) $, hence it allows extensions in $ \mathcal{A}_M(S_\gamma) $ in the following situations:\\
- when $ \gamma(M) $ is an integer and $ \gamma<\gamma(M)-1 $,\\
- when $ \gamma(M) $ is not a integer and $ \gamma<[\gamma(M)] $ (the
brackets denote the integer part).\\ 
This result is not optimal. In particular, the case of all sequences 
$ M $ with $ \gamma(M)\leq 1 $ is not covered. In fact, for such sequences, it is not possible to deduce from \cite{SV} whether the map $ B\, :\, \mathcal{A}_M(S_\gamma)\longrightarrow \Lambda_M(\mathbb{N}) $ is surjective for some $ \gamma>0 $ or not.

However, there are some indications.
For Gevrey sequences, we have recalled the characterization of
surjectivity $ \gamma<\gamma(M) $, and one can even find linear continuous extension operators as mentioned above: this is theorem 5.10 of
\cite{SV}, which is based on a refinement of the classical methods of Laplace transforms. 
Thus, the Gevrey case suggests that for any strongly regular sequence $ M $, one should have corresponding extension operators
as soon as $ \gamma<\gamma(M) $ (and, of course, generally not for any larger $ \gamma $).
This expectation will be satisfied in what follows.

\subsection{The main theorem}\label{myext}

Our approach is based on the following fact: for Gevrey sequences,
a statement quite similar to theorem 5.10 of \cite{SV} had been obtained previously in \cite{VT1}, with a completely different scheme of proof and with an extra assumption, namely $ \gamma<2 $. Here, the flat functions of theorem \ref{sectflat} will allow us to extend the method of \cite{VT1} to general strongly regular sequences $ M $. We shall also overcome the additional limitation on $ \gamma $.
 
\begin{thm}\label{mythm}
Let $ M $ be a strongly regular sequence and let $ \gamma $ be a real number with
$ 0<\gamma<\gamma(M) $. One can then find a constant
$ d\geq 1 $, depending only on $ M $ and $ \gamma $,
such that, for any real $ \sigma>0 $, there exists a linear continuous operator
\begin{equation*}
T_{\gamma,\sigma}\, :\, \Lambda_{M,\sigma}(\mathbb{N})\longrightarrow \mathcal{A}_{M,d\sigma}(S_\gamma)
\end{equation*}
satisfying $ BT_{\gamma,\sigma}\lambda=\lambda $ for any element $ \lambda $ of
$ \Lambda_{M,\sigma}(\mathbb{N}) $.
\end{thm}
\begin{proof} We distinguish two cases in the proof.
\begin{center}
\textsf{first case}: $ \gamma<2 $
\end{center}
In this case, $ S_\gamma $ is a subsector of the complex plane $ \mathbb{C} $ instead of the Riemann surface $ \Sigma $. Putting $ z=x+iy $ for $ z\in\mathbb{C} $, we identify $ \mathbb{C} $ and $ \mathbb{R}^2 $ in the standard way and we denote by
$ \bar\partial $ the Cauchy-Riemann operator
$ \frac{1}{2}\big(\frac{\partial}{\partial x}+
i\frac{\partial}{\partial y}\big) $. Let
$ D $ and $ D' $ be two open discs centered at $ 0 $, with
$ \overline{D}\subset D' $. For a given $ \sigma>0 $, let $ \chi $ 
be a function belonging to $ \mathcal{C}_{M,\sigma}(\mathbb{C}) $, supported in $ D' $ and identically equal to $ 1 $ in $ D $ (for instance, a cut-off function of Bruna type \cite{BBMT}\cite{Bruna}\cite{CC2}). The proof can now be cut into several steps.

\noindent (i) \textit{Construction of ultradifferentiable extensions in $ \mathbb{C} $ with formal
holomorphy at $ 0 $.}
With any given element $ \lambda $ of $ \Lambda_{M,\sigma}(\mathbb{N}) $, we
associate $ \lambda^\mathbb{C}=(\lambda^\mathbb{C}_{jk})_{(j,k)\in\mathbb{N}^2} $ obtained by the natural complexification
\begin{equation}\label{complex}
\sum_{(j,k)\in\mathbb{N}^2}\lambda^\mathbb{C}_{jk}\frac{x^jy^k}{j!k!}=
\sum_{\ell\in\mathbb{N}}\lambda_\ell\frac{(x+iy)^\ell}{\ell!},
\end{equation}
which amounts to putting $ \lambda^\mathbb{C}_{jk}=i^k\lambda_{j+k} $. Remark that the map
$ \lambda\longmapsto \lambda^\mathbb{C} $ acts as a linear continuous operator $ \Lambda_{M,\sigma}(\mathbb{N})\longrightarrow \Lambda_{M,\sigma}(\mathbb{N}^2) $, with
norm $ 1 $.
Then we put  
$ g_\lambda= E_\sigma \lambda^\mathbb{C} $, where 
$ E_\sigma $ is the extension map $ \Lambda_{M,\sigma}(\mathbb{N}^2)
\longrightarrow\mathcal{C}_{M,b\sigma}(\mathbb{C}) $ given by proposition
\ref{extclassic}, and chosen in such a way that the extensions
are supported in $ D $. 
By \eqref{complex}, it is clear that  
$ \bar\partial g_\lambda $ is flat at $ 0 $. Therefore, proceeding as for
\eqref{udflat1}, we get, for any $ K\in\mathbb{N}^2 $ and any
$ z\in\mathbb{C} $, 
\begin{equation}\label{dbarflat}
\big\vert D^K\big(\bar\partial g_\lambda\big)(z)\big\vert\leq c_3\nu_1(\sigma)\vert \lambda\vert_\sigma (c_4\sigma)^k
k!M_k h_M(c_4\sigma\vert z\vert),
\end{equation}
where $ \nu_1(\sigma) $ denotes
the operator norm of $ E_\sigma $, $ c_3=4A^2M_1 $ and $ c_4=4A^2b $.

\noindent (ii) \textit{Division of a flat ultradifferentiable function by a flat ultraholomorphic function.}
For any real $ \tau>0 $ and any $ z\in S_\gamma $, put now $ \psi_\tau(z)=
G(\tau z) $, where $ G $ is the function of theorem
\ref{sectflat} and lemma \ref{flatuplo}. Using \eqref{flatlo}, \eqref{dbarflat} and
lemma \ref{elem}, we get, for
any bi-index $ L\in\mathbb{N}^2 $ and any point $ z $ of $ S_\gamma $, the estimate
\begin{equation*}
\left\vert D^L\Big(\frac{1}{\psi_\tau}\bar\partial g_\lambda\Big)(z)\right\vert\leq b_{15}c_3\nu_1(\sigma)\vert\lambda\vert_\sigma\left(\frac{h_M(c_4\sigma\vert z\vert)}{h_M(b_{17}\tau\vert z\vert)}\right)(b_{16}\tau+c_4\sigma)^\ell \ell! M_\ell.
\end{equation*}
Now we put
\begin{equation}\label{choose}
\tau=\rho(2)c_4b_{17}^{-1}\sigma\quad\text{and}\quad \psi=\psi_\tau.
\end{equation}
Thanks to \eqref{hmodg}, the preceding
estimate then yields
\begin{equation}\label{quotient}
\left\vert D^L\Big(\frac{1}{\psi}\bar\partial g_\lambda\Big)(z)\right\vert\leq \nu_2(\sigma)\vert \lambda\vert_\sigma (c_5\sigma)^\ell\, \ell!M_\ell\, h_M(c_5\sigma\vert z\vert)
\end{equation}
for some suitable positive constants $ \nu_2(\sigma) $ and $ c_5 $.

\noindent (iii) \textit{Solution of a $ \bar\partial $-problem.}
From \eqref{quotient}, we see in particular that $ \frac{1}{\psi}\bar\partial g_\lambda $ belongs to $ \mathcal{C}_{M,c_5\sigma}(S_\gamma) $ and that its norm in this space
is majorized by $ \nu_2(\sigma)\vert \lambda\vert_\sigma $. 
We apply part (ii) of proposition \ref{extclassic}, with $ \Omega=S_\gamma\cap D' $.
Put $ v_\lambda=F_{c_5\sigma}
\big(\frac{1}{\psi}\bar\partial g_\lambda\big) $. Then 
$ v_\lambda $ belongs to $ \mathcal{C}_{M,cc_5\sigma}(\mathbb{C}) $ and we have $ \Vert v_\lambda\Vert_{cc_5\sigma}\leq \nu_3(\sigma)\vert \lambda\vert_\sigma $ where $ \nu_3(\sigma) $ is the product of
$ \nu_2(\sigma) $ and of the operator norm of $ F_{c_5\sigma} $. 
Since $ v_\lambda $ coincides
with $ \frac{1}{\psi}\bar\partial g_\lambda $ on
$ S_\gamma\cap D' $ and $ \frac{1}{\psi}\bar\partial g_\lambda $
vanishes on $ S_\gamma\backslash D $, we have actually
\begin{equation}\label{coinc}
\chi v_\lambda = \frac{1}{\psi}\bar\partial g_\lambda\ \ \textrm{ in all of }\ S_\gamma.
\end{equation}
Moreover, lemma \ref{elem} shows that $ \chi v_\lambda $ 
belongs to $ \mathcal{C}_{M,c_6\sigma}(\mathbb{C}) $ with $ c_6=cc_5+1 $, and that
$ \Vert\chi v_\lambda\Vert_{\mathbb{C},c_6\sigma}\leq 
\nu_4(\sigma)\vert \lambda\vert_\sigma $ for some suitable $ \nu_4(\sigma) $. 
Consider now the convolution $ u_\lambda=\mathcal{K}*(\chi v_\lambda) $, where $ \mathcal{K} $ denotes the Cauchy kernel $ \mathcal{K}(\zeta)=(\pi\zeta)^{-1} $. Since $ \chi v_\lambda $ is
compactly supported in $ D $, the function $ u_\lambda $ solves $ \bar\partial u_\lambda=\chi v_\lambda $ in $ \mathbb{C} $. Moreover, for any $ L\in\mathbb{N}^2 $ it is routine to check that 
$\sup_{z\in\mathbb{C}}\vert D^L u_\lambda(z)\vert\leq \pi\sup_{\zeta\in D}\vert D^L(\chi v_\lambda)(\zeta)\vert$, hence
\begin{equation}\label{estimu} 
\sup_{z\in\mathbb{C}}\vert D^L u_\lambda(z)\vert\leq 
\pi\nu_4(\sigma)\vert\lambda\vert_\sigma (c_6\sigma)^\ell \ell!M_\ell.
\end{equation}

\noindent (iv) \textit{Addition of a flat correction to obtain a holomorphic extension.}
Using \eqref{flatup}, \eqref{choose}, \eqref{estimu} and lemma \ref{elem}, we derive 
that $ \psi u_\lambda $ belongs to $ \mathcal{C}_{M,c_7\sigma}(S_\gamma) $ with $ c_7=\rho(2)c_4b_{13}b_{17}^{-1}+c_6 $, and that we have
$ \Vert \psi u_\lambda\Vert_{S_\gamma,c_7\sigma}\leq 
\pi\nu_4(\sigma)
\vert \lambda\vert_\sigma $.
Put $ f_\lambda= g_\lambda-\psi u_\lambda. $ Then $ f_\lambda $ is
well-defined and holomorphic in $ S_\gamma $ since $ \bar\partial f_\lambda=
\bar\partial g_\lambda-\psi\bar\partial u_\lambda=
\bar\partial g_\lambda-\psi\chi v_\lambda=0 $ in $ S_\gamma $, thanks to \eqref{coinc}.
A quick look at the previous constructions also shows that the functions  $ g_\lambda $, $ v_\lambda $, $ u_\lambda $, and subsequently $ f_\lambda $, all depend linearly on $ \lambda $. Our estimates on $ g_\lambda $ and $ \psi u_\lambda $ show moreover that the
map $ \lambda\longmapsto f_\lambda $ is continuous from 
$ \Lambda_{M,\sigma}(\mathbb{N}) $ to
$ \mathcal{A}_{M,c_8\sigma}(S_\gamma) $ where $ c_8=\max(b,c_7) $ depends only on
$ \gamma $ and on the sequence $ M $.
Finally, since $ \psi $, and consequently $ \psi u_\lambda $, are flat at the origin, we have $ f_\lambda^{(j)}(0)=(\partial^j f_\lambda/\partial x^j)(0)=
(\partial^j g_\lambda/\partial x^j)(0)=\lambda_j $ for any integer $ j\geq 0 $. Thus, it suffices to put $ T_{\gamma,\sigma}\lambda=f_\lambda $ and $ d=c_8 $ to get the desired conclusion in this case. 
\begin{center}
\textsf{second case}: $ \gamma\geq 2 $
\end{center}
We use the naive idea of reducing this case to the first one by means of a suitable ramification. The less obvious part consists in showing that the required estimates are preserved by this process. Pick an integer $ q $ such that $ \gamma/q <2 $ (hence $ q\geq 2 $). We shall use the estimate
\begin{equation}\label{Mqj}
M_j\leq (M_{qj})^{1/q}\leq A^jM_j\text{ for any }j\in\mathbb{N},
\end{equation}
which is immediate by \eqref{Mmodg} and \eqref{Mfast}. Now, with any sequence
$ \lambda $ in
$ \Lambda_{M,\sigma}(\mathbb{N}) $, we associate the sequence $ \lambda^* $ given by
\begin{equation*}
\lambda^*_{qj}=\lambda_j\frac{(qj)!}{j!}\quad\text{and}\quad\lambda^*_{qj+k}=0
\quad\text{for any }j\in\mathbb{N}\text{ and any }k=1,\ldots,q-1.
\end{equation*}
By \eqref{Mqj}, it is easy to see that $ \lambda^* $ belongs to 
$ \Lambda_{M^{1/q},\sigma^{1/q}}(\mathbb{N}) $ and that its norm in this space is majorized by $ \vert\lambda\vert_\sigma $. Recall also that we have
$ \gamma/q<\gamma(M)/q=\gamma(M^{1/q}) $ by \eqref{gamma}. We can therefore extend $ \lambda^* $ by applying the first case of the proof with $ M $ replaced by $ M^{1/q} $ and $ \gamma $ replaced by $ \gamma/q $. We obtain thus a function $ h_\lambda $ 
in $ \mathcal{A}_{M^{1/q},c_9\sigma^{1/q}}(S_{\gamma/q}) $ depending 
linearly and continuously on $ \lambda \in\Lambda_{M,\sigma}(\mathbb{N}) $ and
verifying
\begin{equation}\label{derivg}
Bh_\lambda=\lambda^*.
\end{equation}
Of course, the constant $ c_9 $ depends only on $ M $ and $ \gamma $.
For $ z\in S_\gamma $, we put now $ f_\lambda(z)=h_\lambda(z^{1/q}) $. Obviously, the function $ f_\lambda $ is holomorphic in $ S_\gamma $ and bounded by the supremum norm of
$ h_\lambda $. In order to estimate its derivatives, we shall follow the general pattern  of \cite{VT2}. Consider the differential operator
$ \mathbf{Y} = q^{-1}w^{1-q}\frac{\partial}{\partial w} $ on $ \mathbb{C}\setminus\{0\} $, so that
\begin{equation}\label{push}
f_\lambda^{(\ell)}(w^q)=(\mathbf{Y}^\ell h_\lambda)(w)\text{ for any }w\in S_{\gamma/q}
\text{ and any integer }\ell\geq 1.
\end{equation}
Proceeding by induction on $ \ell $ as in the proof of proposition 2.5 of \cite{VT2} (but in a much simpler situation), we obtain
\begin{equation*}
\mathbf{Y}^\ell=\sum_{k=1}^\ell\mathrm{Y}_{\ell,k}(w)
\frac{\partial^k}{\partial w^k}
\end{equation*}
with
\begin{equation}\label{iterates}
\left\vert\frac{\partial^j \mathrm{Y}_{\ell,k}}{\partial w^j}(w)
\right\vert\leq
\left(4q^{-1}\right)^\ell 2^{\ell+j-k}(\ell+j-k)!\, \vert w\vert^{k-j-q\ell}
\end{equation}
for any integers $ j $, $k $, $ \ell $ with $ 1\leq k\leq\ell $ and $ j\geq 0 $. Now we introduce the polynomial
\begin{equation}\label{defP}
\mathcal{P}_{\lambda,\ell}(w)=\sum_{j=0}^{\ell-1}\lambda_j\frac{w^{qj}}{j!}.
\end{equation}
From \eqref{derivg} and from the definition of $ \lambda^* $, we see that
\begin{equation*}
\mathcal{P}_{\lambda,\ell}(w)=\sum_{p=0}^{q\ell-1}h_\lambda^{(p)}(0)\frac{w^p}{p!}.
\end{equation*}
The Taylor formula for $ h_\lambda $ between $ 0 $ and any point $ w $ in
$ S_{\gamma/q} $ yields therefore
\begin{equation*}
\begin{split}
\left\vert\frac{\partial^k}{\partial w^k}\big(h_\lambda-\mathcal{P}_{\lambda,\ell}\big)(w)\right\vert & \leq\sup_{\zeta\in ]0,w[}\big\vert h_\lambda^{(q\ell)}(\zeta)\big\vert\frac{\vert w\vert^{q\ell-k}}{(q\ell-k)!} \\
 & \leq \nu_5(\sigma)\vert\lambda\vert_\sigma
(c_9\sigma^{1/q})^{q\ell}\frac{(q\ell)!}{(q\ell-k)!}(M_{q\ell})^{1/q}\vert w\vert^{q\ell-k},
\end{split}
\end{equation*}
where $ \nu_5(\sigma) $ denotes the operator norm of the map $ \lambda\longmapsto h_\lambda $ from $ \Lambda_{M,\sigma}(\mathbb{N}) $ to $ \mathcal{A}_{M^{1/q},c_9\sigma^{1/q}}(S_{\gamma/q}) $. Using \eqref{Mqj} and the elementary estimate
$ (q\ell)!\leq 2^{q\ell}(q\ell-k)!k! $, we derive
\begin{equation*}
\left\vert\frac{\partial^k}{\partial w^k}\big(h_\lambda-\mathcal{P}_{\lambda,\ell}\big)(w)\right\vert \leq \nu_5(\sigma)\vert\lambda\vert_\sigma
(c_{10}\sigma)^\ell\, k!\, M_\ell\, \vert w\vert^{q\ell-k},
\end{equation*}
with $ c_{10}=(2c_9)^qA $. Together with \eqref{iterates} and the obvious fact $ (\ell-k)!k!\leq \ell! $, this yields finally
\begin{equation}\label{piece1}
\big\vert\big(\mathbf{Y}^\ell(h_\lambda-\mathcal{P}_{\lambda,\ell}\big)(w)\big\vert\leq \nu_5(\sigma)\vert\lambda\vert_\sigma
(c_{11}\sigma)^\ell \ell! M_\ell\quad\text{for any } w\in S_{\gamma/q},
\end{equation}
with $ c_{11}=8c_{10}/q $. Now recall from \eqref{defP} that $ \mathcal{P}_{\lambda,\ell}(w) $ can be written as $ \mathcal{Q}(w^q) $ where $ \mathcal{Q} $ is a polynomial of degree at most $ \ell-1 $, hence
\begin{equation}\label{piece2}
\big(\mathbf{Y}^\ell\mathcal{P}_{\lambda,\ell}\big)(w)=\mathcal{Q}^{(\ell)}(w^q)=0.
\end{equation}
Gathering \eqref{push}, \eqref{piece1} and \eqref{piece2}, we see that $ f_\lambda $ belongs to $ \mathcal{A}_{M,d\sigma}(S_\gamma) $ with $ d=c_{11} $, and that it depends linearly and continuously on
$ \lambda $.  At last, we know from \eqref{derivg} that the Taylor series of $ h_\lambda $ at $ 0 $ is 
$ \sum_{j\in\mathbb{N}}\lambda_j\frac{w^{qj}}{j!} $, which means that the expansion of $ f_\lambda $ is given by $ \sum_{j\in\mathbb{N}}\lambda_j\frac{w^j}{j!} $. Therefore we have $ Bf_\lambda=\lambda $ and the proof is complete. 
\end{proof}

\subsection{Comments on the strong regularity assumption}\label{gene}

The problem of sectorial extensions could be put 
under weaker assumptions on $ M $, say \eqref{Mnorm} and
\eqref{Mlogc}. But as
in the $ \mathcal{C}^\infty $ case, it is easy to see
that the desired surjectivity property implies also the
surjectivity of $ B\, :\, \mathcal{C}_M(\mathbb{R})\longrightarrow \Lambda_M(\mathbb{N}) $. Therefore, Petzsche's results \cite{Pet} show that assumption \eqref{Msnqa} is necessary. 
The additional moderate growth assumption \eqref{Mmodg} seems more related to technical reasons. Consider, for example, $ M_j=e^{j^2} $. In this case, for which \eqref{Mnorm},
\eqref{Mlogc} and \eqref{Msnqa} hold, but not \eqref{Mmodg}, it is easy to
check that theorem 5.6 of \cite{SV} applies: for any real $ \gamma>0 $, there
exists a continuous extension map $ T_\gamma $ from the (LB)-space 
$ \Lambda_M(\mathbb{N}) $ to the
(LB)-space $ \mathcal{A}_M(S_\gamma) $. The arbitrary aperture agrees
with the fact that $ (\mathsf{P}_\gamma) $ holds here for any
$ \gamma $.

In this particular example, disregarding continuity properties, one can also obtain an extension procedure working simultaneously for every $ \gamma $. Indeed, being
given a sequence $ \lambda$ in $ \Lambda_M(\mathbb{N}) $, section 2 of \cite{Zhang}
provides a function holomorphic
in a whole ``punctured disc'' $ D=\{z\in\Sigma\, ;\, \vert z\vert<\delta \} $ on the
Riemann surface $ \Sigma $, and whose restriction to every bounded sector
$ S_\gamma\cap D $ belongs to $ \mathcal{A}_M(S_\gamma\cap D) $ and
satisfies $ Bf=\lambda $. We do not know
whether a similar statement holds for more general sequences $ M $ satisfying \eqref{Mnorm}, \eqref{Mlogc} and \eqref{Msnqa}, but not \eqref{Mmodg}.

\subsection{The case of Beurling classes}\label{Beurling}

Instead of a Carleman class $ \mathcal{A}_M(S_\gamma) $, it is also possible to consider a  Beurling class, that is
the Fr\'echet space $ \mathcal{A}^-_M(S_\gamma) $ obtained as the
projective limit of spaces $ \mathcal{A}_{M,\sigma}(S_\gamma) $. In the same
way, one defines a Fr\'echet
space $ \Lambda^-_M(\mathbb{N}) $ as the projective limit of spaces $ \Lambda_{M,\sigma}(\mathbb{N}) $ and we have an induced Borel map
$ B\, :\, \mathcal{A}^-_M(S_\gamma)\longrightarrow\Lambda^-_M(\mathbb{N}) $. 
Theorem \ref{mythm} then has the following corollary.

\begin{cor}
For any strongly regular sequence $ M $ and any real number $ \gamma $ with
$ 0<\gamma<\gamma(M) $, the Borel map
$ B\, :\, \mathcal{A}^-_M(S_\gamma)\longrightarrow\Lambda^-_M(\mathbb{N}) $ is
surjective.
\end{cor}
\begin{proof} The result is an easy consequence of theorem \ref{mythm} and of the following argument, inspired by \cite{CC2}: for any element
$ \lambda $ of $ \Lambda^-_M(\mathbb{N}) $, one can find a strongly regular sequence
$ N $, with $ \gamma<\gamma(N) $, such that the sequence $ \lambda $ belongs to 
$ \Lambda_N(\mathbb{N}) $ and the Carleman class $ \mathcal{A}_N(S_\gamma) $ is contained in the Beurling class $ \mathcal{A}^-_M(S_\gamma) $. The construction of $ N $ is a slight variation on lemma 16 and proposition 17 of \cite{CC2}; we will not describe the details here. 
\end{proof} 

It should be emphasized that Schmets and Valdivia have also considered in \cite{SV} the case of Beurling classes, quite in the same spirit as for Carleman classes. A nice feature of theorem 4.5 of \cite{SV} is that it provides linear continuous extension maps.

\section{Division by flat functions}\label{div}

\subsection{Setting of the problem}\label{divback}

Denote by $ \mathcal{C}^\infty(\mathbb{R}^n,0) $ the ring of $ \mathcal{C}^\infty $ function germs
at the origin of $ \mathbb{R}^n $. For any open neighborhood $ \Omega $ of the origin,  denote by $ \pi_\Omega $ the canonical mapping which, to each element of
$ \mathcal{C}^\infty(\Omega) $, associates its germ in $ \mathcal{C}^\infty(\mathbb{R}^n,0) $. Let $ \mathcal{I} $ be an ideal of $ \mathcal{C}^\infty(\mathbb{R}^n,0) $. An ideal $ \mathcal{I}_\Omega$ of $ \mathcal{C}^\infty(\Omega) $ is called a \emph{representative} of $ \mathcal{I} $ if it satisfies $ \pi_\Omega(\mathcal{I}_\Omega)=\mathcal{I} $, and we say that $ \mathcal{I} $ is \emph{closed} if, for any sufficiently small $ \Omega $, it has a closed representative in the Fr\'echet space $ \mathcal{C}^\infty(\Omega) $.

Now let $ X $ be a germ of closed subset at the origin of $ \mathbb{R}^n $. With the usual confusion between germs and their representatives, an element of $ \mathcal{C}^\infty(\mathbb{R}^n,0) $ is said to be \emph{flat on $ X $} if it vanishes, together with all its derivatives, on $ X $. We denote by $ \imax^\infty_X $ the ideal of all such germs. A classical result of Tougeron (\cite{Tougbook}, proposition V.2.3), when stated from a local viewpoint, asserts that
for any closed ideal $ \mathcal{I} $ of $ \mathcal{C}^\infty(\mathbb{R}^n,0) $, 
the ideal of elements of $ \mathcal{I} $ which
are flat on $ X $ is generated over $ \mathcal{I} $ by $ \imax^\infty_X$, which can be written 
\begin{equation}\label{eqtoug}
\mathcal{I}\cap \imax^\infty_X=\imax^\infty_X\mathcal{I}.
\end{equation}

In particular, one has $ \imax^\infty_X=\imax^\infty_X\imax^\infty_X $. It should be
emphasized that these properties are delicate even when $ X=\{0\} $.
Interesting in themselves, they also have a number of applications in differential analysis. This motivates their
study in the setting of ultradifferentiable classes. 

\subsection{The ultradifferentiable case}

Let $ M $ be a strongly regular sequence. We consider the ring
$ \mathcal{C}_M(\mathbb{R}^n,0) $ of all those germs $ f $ of
$ \mathcal{C}^\infty(\mathbb{R}^n,0) $ which have a representative in 
$ \mathcal{C}_M(\Omega_f) $ for some open neighborhood $ \Omega_f $ of
$ 0 $. Just as in the $ \mathcal{C}^\infty $ case, being given an open neighborhood $ \Omega $ of $ 0 $ and an
ideal $ \mathcal{I}_M $ of $ \mathcal{C}_M(\mathbb{R}^n,0) $, we say that an ideal $ \mathcal{I}_{M,\Omega} $ of $ \mathcal{C}_M(\Omega) $ is a \emph{representative} of $ \mathcal{I}_M $ if it satisfies $ \pi_\Omega(\mathcal{I}_{M,\Omega})=\mathcal{I}_M $, and we say that $ \mathcal{I}_M $ is \emph{closed} if, for any sufficiently small $ \Omega $, it has a closed representative in the 
(LB)-space $ \mathcal{C}_M(\Omega) $. Finally, we denote by $ \imax^\infty_{X,M} $ the ideal of germs of
$ \mathcal{C}_M(\mathbb{R}^n,0) $ which are flat on the germ of closed subset $ X $.

\begin{lem}\label{flatuplo2}
Let $ M $ be a strongly regular sequence and let $ V $ be a vector subspace of
$ \mathbb{R}^n $.
Then, for any real $ \tau>0 $, there exists 
a real non-negative function
$ v_\tau $ which belongs to $ \imax^\infty_{V,M}(\mathbb{R}^n) $, has
no zero in $ \mathbb{R}^n\setminus V $ and satisfies, for any multi-index
$ J $ and any point $ x $ in $ \mathbb{R}^n\setminus V $, 
\begin{equation}\label{flatlo2}
\Big\vert D^J\Big(\frac{1}{v_\tau}\Big)(x)\Big\vert \leq d_1(d_2\tau)^j j!M_j
\big(h_M\big(d_3\tau \dist(x,V)\big)\big)^{-1}
\end{equation}
with positive constants $ d_1 $, $ d_2 $, $ d_3 $ depending only on $ M $ and $ X $.
\end{lem}
\begin{proof} Using the natural identification
of $ \mathbb{R}^n $ as a totally real subset
of $ \mathbb{C}^n $, we define, for any real $ \varepsilon>0 $, the set
$ V_\varepsilon=\{\zeta\in\mathbb{C}^n\, ;\, \vert\Im\zeta\vert<\varepsilon \dist(\zeta,V)\} $. Note that $ V_\varepsilon $ contains $ \mathbb{R}^n\setminus V $.
After a suitable linear change of coordinates, we can
assume $ V=W\cap\mathbb{R}^n $ with 
$ W=\{\zeta\in\mathbb{C}^n\, ;\, \zeta_1=\cdots=\zeta_k=0\} $ and $ k=\codim V $. For
$ \zeta\in\mathbb{C}^n $, put
$ \xi=\Re\zeta $. Then 
\begin{equation}\label{ex1}
\dist(\zeta,V)^2=\xi_1^2+\cdots+\xi_k^2+\vert\Im\zeta\vert^2. 
\end{equation}
Consider $ Q(\zeta)=\zeta_1^2+\cdots+\zeta_k^2 $. For $ 1\leq j\leq k $, one
has $ \vert\xi_j\vert\leq\vert\zeta_j\vert\leq \dist(\zeta,V) $ by \eqref{ex1}. Assuming
now that
$ \zeta $ belongs to $ V_\varepsilon $, one has also $ \vert\zeta_j-\xi_j\vert<\varepsilon \dist(\zeta,V) $ and we obtain $ \vert Q(\zeta)-Q(\xi)\vert< 2k\varepsilon \dist(\zeta,V)^2 $. Using \eqref{ex1} again, one gets
$ Q(\xi)\geq (1-\varepsilon^2)\dist(\zeta,V)^2 $. All this yields
\begin{equation}\label{realpart}
\Re Q(\zeta)>(1-\varepsilon^2-2k\varepsilon)\dist(\zeta,V)^2.
\end{equation}
Now let $ \gamma $ be a real number with $ 0<\gamma<\gamma(M) $.
Pick $ \delta $ with $ 0<\delta<\min(1,2\gamma) $ and choose $ \varepsilon>0 $ 
small enough
to have $ 1-\varepsilon^2-2k\varepsilon>\cos(\delta\frac{\pi}{2}) $.
Denote by $ \sqrt{\, } $ the natural
determination of the square root in $ \mathbb{C}\setminus ]-\infty,0] $.
Since $ \cos(\delta\frac{\pi}{2})>0 $ and $ \dist(\zeta,V)^2\geq \vert\zeta\vert^2\geq\vert Q(\zeta)\vert $, the estimate \eqref{realpart} shows that
the function $ \Phi $ given by $ \Phi(\zeta)=\sqrt{Q(\zeta)} $
is well-defined and holomorphic in $ V_\varepsilon $ and that it satisfies $ \Phi(V_\varepsilon)\subset S_{\delta/2} $, hence
\begin{equation}\label{phi1}
\Phi(V_\varepsilon)\subset S_\gamma.
\end{equation}
One has also clearly, for some suitable constant $ d_4\geq 1 $, 
\begin{equation}\label{phi2}
d_4^{-1}\dist(\zeta,V)\leq \vert\Phi(\zeta)\vert\leq d_4\dist(\zeta,V)\text{ for any }\zeta\in V_\varepsilon.
\end{equation}
For $ \zeta\in V_\varepsilon $, we consider 
$ H(\zeta)=G(\tau\Phi(\zeta)) $, where $ G $ denotes the function
of theorem \ref{sectflat}. We shall see that the function
$ v_\tau $ defined by
\begin{equation}\label{v1}
v_\tau(x)=H(x) \text{ for } x\in\mathbb{R}^n\setminus V,\quad v_\tau(x)=0
\text{ for } x\in V
\end{equation}
has all the required properties. For any $ x $ in $ \mathbb{R}^n\setminus V $, consider the
closed polydisc $ P_x=\{\zeta\in\mathbb{C}^n\, ;\, \vert \zeta_j-x_j\vert\leq
\frac{\varepsilon}{2\sqrt{n}}\dist(x,V)\text{ for }j=1,\dots,n\}. $ Then $ P_x $ is contained in $ V_\varepsilon $. One can assume
$ \varepsilon<1 $, so that
$ \frac{1}{2}\dist(x,V)\leq \dist(\zeta,V)\leq 2\dist(x,V) $ for any $ \zeta\in P_x $. Since $ H $ is holomorphic in $ V_\varepsilon $, it
restricts to a $ \mathcal{C}^\infty $ function in $ \mathbb{R}^n\setminus V $ and its derivatives at the point $ x $ can be estimated by the Cauchy formula on $ P_x $. The scheme of proof is then essentially
the same as in lemma \ref{flatuplo}: taking into account the upper bounds in \eqref{sharp} and \eqref{phi2}, we get the estimate
$ \vert D^JH(x)\vert\leq d_5^j\, j!\, \dist(x,V)^{-j}
h_M\big(d_6\tau \dist(x,V)\big) $ for any multi-index $ J $, 
with $ d_5=2\sqrt{n}\varepsilon^{-1} $ and $ d_6=2\kappa_3d_4 $. Using \eqref{hmodg}, we derive 
\begin{equation}\label{v2}
\vert D^JH(x)\vert\leq (d_7\tau)^j j!M_j
h_M\big(d_8\tau \dist(x,V)\big)\text{ for any }J\in\mathbb{N}^n,
\end{equation}
with $ d_7=d_5d_6\rho(2)\varepsilon^{-1} $ and $ d_8=d_6\rho(2) $. In particular, we see
that all the derivatives $ D^JH(x) $ tend to $ 0 $ as $ x $ approaches $ V $ in $ \mathbb{R}^n\setminus V $.
By Hestenes lemma, the function $ v_\tau $ defined in \eqref{v1} is therefore
$ \mathcal{C}^\infty $ in $ \mathbb{R}^n $, and by \eqref{v2} it belongs to
$ \imax^\infty_{V,M}(\mathbb{R}^n) $ as announced.
The derivation of \eqref{flatlo2} goes along the same lines: using the lower bounds in
\eqref{sharp} and \eqref{phi2}, the Cauchy formula on $ P_x $ yields the
desired estimate \eqref{flatlo2} by virtue of \eqref{hmodg}.
\end{proof}

We can now state the key result of this section.

\begin{prop}\label{divis}
Let $ X $ be a germ of real-analytic submanifold at the
origin in $ \mathbb{R}^n $. 
Then for any strongly regular sequence $ M $ and any finite family
$ u_1,\dots,u_p $ of germs belonging to $ \imax^\infty_{X,M} $, one can find
an element $ v $ of $ \imax^\infty_{X,M} $ whose germ of zero set is precisely
$ X $ and such that
$ u_i $ belongs to $ v\,\imax^\infty_{X,M} $ for $ i=1,\dots,p $.
\end{prop}
\begin{proof} After a suitable real-analytic change of coordinates in a
neighborhood $ \Omega $ of $ 0 $, one can assume that $ X=\Omega\cap V $, where
$ V $ is a vector subspace of $ \mathbb{R}^n $. 
The proof consists in showing, from lemma \ref{flatuplo2} and from the flatness
of the $ u_i $, that it is possible to choose $ v $ as the germ of $ v_\tau $ for some
suitable $ \tau $. We shall not describe all the details, since the arguments essentially mimic step (ii) of the proof of theorem \ref{mythm}. First, the flatness of the $ u_i $ yields, for any multi-index $ K\in\mathbb{N}^n $, any $ x $ in $ \Omega $ and any $ i=1,\ldots,p $,
\begin{equation}\label{uiflat}
\vert D^Ku_i(x)\vert\leq d_9d_{10}^k k!M_kh_M\big(d_{11}\dist(x,X)\big)
\end{equation}
for some suitable positive constants $ d_9 $, $ d_{10} $, $ d_{11} $ (the proof is the same as for \eqref{udflat1}, except that the Taylor formula
is used between $ x $ and a point $ \hat{x} $ in $ X $ satisfying $ \dist(x,X)=\vert x-\hat{x}\vert $). Without loss of generality, one can assume that $ \Omega $ and $ V $ intersect transversally, so that any point $ x $ of $ \Omega $ verifies $ \dist(x,X)\leq d_{12}\dist(x,V) $ for some positive constant $ d_{12} $ depending only on $ \Omega $ and $ V $. Hence \eqref{flatlo2} yields 
\begin{equation}\label{flatlo3}
\Big\vert D^J\Big(\frac{1}{v_\tau}\Big)(x)\Big\vert \leq d_2(d_3\tau)^j j!M_j
\big(h_M\big(d_{13}\tau \dist(x,X)\big)\big)^{-1}
\end{equation} 
for any multi-index $ J $ in $ \mathbb{N}^n $ and any $ x $ in $ \Omega $, with $ d_{13}=d_3d_{12}^{-1} $. Now it is enough to choose $ \tau\geq\rho(2)d_{11}d_{13}^{-1} $: putting $ v(x)=v_\tau(x) $ for any $ x $ in
$ \Omega $, the desired result then follows from \eqref{hmodg}, \eqref{uiflat}, \eqref{flatlo3} and lemma \ref{elem}.
\end{proof}

\begin{rem}
Contrarily to what happens in the $ C^\infty $ setting (\cite{Tougbook}, lemma V.2.4), it is easy to see that proposition \ref{divis} is no
longer true if $ u_1,\ldots,u_p $ is replaced by a countable family $ (u_i)_{i\geq 1} $, even in the simplest case $ X=\{0\} $. 
\end{rem}

We obtain finally a result in the spirit of property \eqref{eqtoug}.

\begin{thm}\label{divflat}
Let $ X $ be a germ of real-analytic submanifold at the origin in $ \mathbb{R}^n $. Then, for any  strongly regular sequence $ M $ and 
any closed ideal $ \mathcal{I}_M $ of $ \mathcal{C}_M(\mathbb{R}^n,0) $, one has
\begin{equation*}
\mathcal{I}_M\cap \imax^\infty_{X,M}=\imax^\infty_{X,M}\mathcal{I}_M.
\end{equation*} 
In particular, $ \imax^\infty_{X,M}=\imax^\infty_{X,M}\imax^\infty_{X,M} $.
\end{thm}
\begin{proof}
It suffices to copy the proof of \eqref{eqtoug} in \cite{Tougbook}, using 
proposition \ref{divis} (only the case $ p=1 $ is required) instead of lemma V.2.4 of \cite{Tougbook}, and the $ \mathcal{C}_M $ version of Whitney's spectral theorem due to Chaumat-Chollet \cite{CC3} instead of the usual $ \mathcal{C}^\infty $ one.
\end{proof}

\begin{prob} 
We do not know whether the previous results are true for more general classes of closed subsets $ X $. In particular, does the identity
$ \imax^\infty_{X,M}=\imax^\infty_{X,M}\imax^\infty_{X,M} $ still hold when $ X $ is a singular real-analytic variety ?
\end{prob}

\begin{flushleft}
Math\'ematiques - B\^atiment M2\\
Universit\'e des Sciences et Technologies de Lille\\
F-59655 Villeneuve d'Ascq Cedex, France\\
e-mail: \texttt{thilliez@agat.univ-lille1.fr}
\end{flushleft}


\begin{thebibliography}{30}

\bibitem{BBMT} \textsc{J. Bonet, R.W. Braun, R. Meise \& B.A. Taylor}, \textit{Whitney's extension theorem for nonquasianalytic classes of ultradifferentiable functions}, Studia Math. \textbf{99} (1991), 155--184.

\bibitem{BrakImm} \textsc{B.L.J. Braaksma \& G.K. Immink}, \textit{A Borel-Ritt theorem with prescribed error bounds}, in: \'Equations diff\'erentielles dans le champ complexe, Vol. I (Strasbourg 1985), Publ. Inst. Rech. Math. Av., Univ. Louis Pasteur, Strasbourg (1988), 1--33.

\bibitem{Bruna} \textsc{J. Bruna}, \textit{An extension theorem of Whitney type for non-quasianalytic classes of functions}, J. London Math. Soc. \textbf{22} (1980), 495--505.

\bibitem{CC0} \textsc{J. Chaumat \& A.-M. Chollet}, \textit{Caract\'erisation et propri\'et\'es des ensembles localement pics de $ A^\infty(D) $}, Duke Math. J. \textbf{47} (1980), 763--787.

\bibitem{CC1} \textsc{J. Chaumat \& A.-M. Chollet}, \textit{Classes de Gevrey non isotropes et application \`a l'interpolation}, Ann. Scuola Norm. Sup. Pisa (IV) \textbf{15} (1988), 615--676.

\bibitem{CC2} \textsc{J. Chaumat \& A.-M. Chollet}, \textit{Surjectivit\'e de l'application restriction \`a un compact dans des classes de fonctions ultradiff\'erentiables}, Math. Ann. \textbf{298} (1994), 7--40.

\bibitem{CC3} \textsc{J. Chaumat \& A.-M. Chollet}, \textit{Caract\'erisation des anneaux noeth\'eriens
de s\'eries formelles \`a croissance contr\^ol\'ee. Application \`a la synth\`ese
spectrale}, Publ. Mat. \textbf{41} (1997), 545--561.

\bibitem{Duren} \textsc{P.L. Duren}, \textit{Theory of $ H^p $ spaces}, Pure and
Applied Mathematics, vol. 38, Academic Press (1970).

\bibitem{Ecalle} \textsc{J. \'Ecalle}, \textit{Les fonctions r\'esurgentes, tome III},
Publications Math\'ematiques d'Orsay, vol. 85-5, Universit\'e de Paris-Sud (1985).

\bibitem{Hof} \textsc{K.J. Hoffman}, \textit{Banach Spaces of Analytic Functions}, Prentice-Hall (1962).

\bibitem{Kom} \textsc{H. Komatsu}, \textit{Ultradistributions, I. Structure theorems and a characterization}, J. Fac. Sci. Tokyo, Sect. IA Math. \textbf{20} (1973), 25--105.

\bibitem{Kor} \textsc{B.I. Korenbljum}, \textit{Conditions of nontriviality of certain
classes of functions analytic in a sector and problems of quasianalyticity}, Soviet Math. Dokl. \textbf{7} (1966), 232--236.

\bibitem{Pet} \textsc{H.-J. Petzsche}, \textit{On E. Borel's theorem}, Math. Ann.  \textbf{282} (1988), 299--313.

\bibitem{Ra} \textsc{J.-P. Ramis}, \textit{D\'evissage Gevrey}, Ast\'erisque \textbf{59-60} (1978), 173--204.

\bibitem{SV} \textsc{J. Schmets \& M. Valdivia}, \textit{Extension maps in ultradifferentiable and ultraholomorphic function spaces}, Studia Math. \textbf{143} (2000), 221--250.

\bibitem{VT1} \textsc{V. Thilliez}, \textit{Extension Gevrey et rigidit\'e dans un secteur}, Studia Math. \textbf{117} (1995), 29--41.

\bibitem{VT2} \textsc{V. Thilliez}, \textit{Sur les fonctions compos\'ees ultradiff\'erentiables}, J. Math. Pures et Appl. \textbf{76} (1997), 499--524.

\bibitem{Tougbook} \textsc{J.-C. Tougeron}, \textit{Id\'eaux de fonctions diff\'erentiables}, Springer Verlag (1972).

\bibitem{Touglectures} \textsc{J.-C. Tougeron}, \textit{An introduction to the theory of Gevrey expansions and to the Borel-Laplace transform with some applications}, Lecture Notes, University of Toronto (1989).

\bibitem{Zhang} \textsc{C. Zhang}, \textit{D\'eveloppements asymptotiques $ q $-Gevrey et s\'eries $ Gq $-sommables}, Ann. Inst. Fourier \textbf{49} (1999), 227--261.

\end{thebibliography}
\end{document}